\documentclass[11pt]{article}

\date{}

\usepackage{epsfig, color}
\usepackage{color,graphicx}
\usepackage{amsfonts}
\usepackage{amsmath}

\newcommand{\bey}{\begin{eqnarray}}
\newcommand{\eey}{\end{eqnarray}}

\newcommand{\beq}{\begin{equation}}
\newcommand{\eeq}{\end{equation}}

\newtheorem{co}{Corollary}[section]

\newtheorem{exam}{Example}[section]

\newcommand{\proof}{\noindent {Proof.} \hspace{2mm}}
\newcommand{\eproof}{$\quad \Box$}
\newtheorem{theorem}{Theorem}[section]

\newtheorem{lemma}{Lemma}[section]
\newtheorem{definition}{Definition}[section]

\newtheorem{remark}{Remark}[section]

\numberwithin{equation}{section}

\begin{document}
\title{A Monotone Finite Volume Method for Time Fractional
Fokker-Planck Equations} \vskip 4mm

\author {Yingjun Jiang\thanks{Department of Mathematics and Scientific
Computing, Changsha University of Science and Technology, Changsha,
410076, China (jiangyingjun@csust.edu.cn).} and Xuejun
Xu\thanks{Institute of Computational Mathematics and
Scientific/Engineering Computing, Academy of Mathematics and Systems
Science, Chinese Academy of Sciences, P.O. Box 2719, Beijing,
100190, P.R. China, and School of Mathematical Sciences, Tongji University, Shanghai, China({\tt xxj@lsec.cc.ac.cn}).}}

\maketitle

\begin{abstract}
 We develop a monotone finite volume method for the time fractional Fokker-Planck equations and theoretically prove its unconditional stability.
 We show that the convergence rate of this method is order 1 in space and if the space grid becomes sufficiently fine, the convergence rate can be improved to order 2. Numerical results are given to support our theoretical findings.   One characteristic of our method is that it has monotone property such that it keeps the
nonnegativity of some physical variables such as density, concentration, etc.

\end{abstract}

\vspace{1.5cm}

\noindent {\bf 2000 Mathematics subject classification:}
 {65M12, 65M06, 35S10}

\vspace{0.5cm}

\noindent {\bf Keywords:} {time fractional
Fokker-Planck equations, finite volume methods}

\section{Introduction}
This paper considers the time fractional Fokker-Planck equation (FFPE)
\begin{equation}\label{model}
\frac{\partial w}{\partial t}=\,_0D^{1-\alpha}_t\left[k_\alpha\frac{\partial^2}{\partial x^2}-\frac{\partial}{\partial x}f(x)\right]w(x,t), \quad a\leq x \leq b,  0<t\leq T,
\end{equation}
subject to the initial condition
\begin{equation}\label{p_b1}
w(x,0)=\varphi(x),\quad a\leq x\leq b,
\end{equation}
and boundary conditions
\begin{equation}\label{p_b2}w(a,t)=g_1(t),\quad w(b,t)=g_2(t),\quad 0<t\leq T,
\end{equation}
where $\alpha\in (0,1)$,  $k_\alpha$ is a positive constant, $f(x),\phi(x),g_1(t),g_2(t)$ are given functions, and $\,_0D^{1-\alpha}_t w(x,t)=\frac1{\Gamma(\alpha)} \frac{\partial}{\partial t}\int_0^t \frac{w(x,s)}{(t-s)^{1-\alpha}}ds$
 is the Riemann-Liouville fractional derivative of order $1-\alpha$ with $\Gamma(x)$ being the gamma function.

 When the equation (\ref{model}) is used to model the sub-diffusion processes in an external force field (see e.g. \cite{Metzler0,Metzler1,So2004,Sokolov2001}),  $k_\alpha$ denotes the generalized diffusion coefficient, and $f(x)=\frac{v'(x)}{m\eta_\alpha}$ represents the external force field governed by the potential $v'(x)$ with $m$ being the mass of the diffusion particle
and $\eta_\alpha$ the generalized friction coefficient.

Under the condition that problem \{(\ref{model})-(\ref{p_b2})\} admits a unique solution $w(x,t)\in C^{2,1}_{x,t}[a,b]\times [0,T]$, the equation (\ref{model}) may be rewritten in the following equivalent form (see e.g. \cite{Chen}):
\begin{equation}\label{modell}
\frac{\partial^\alpha w}{\partial t^\alpha}=\left[k_\alpha\frac{\partial^2}{\partial x^2}-\frac{\partial}{\partial x}f(x)\right]w(x,t), \quad a\leq x \leq b,  0<t\leq T,
\end{equation}
where ${\partial^\alpha w}/{\partial t^\alpha}$ denotes the Caputo fractional derivative of order $\alpha (0<\alpha<1)$ defined by
$$\frac{\partial^\alpha w(x,t)}{\partial t^\alpha}=\frac1{\Gamma(1-\alpha)}\int_0^t\frac{\partial w(x,\eta)}{\partial t}\frac{d\eta}{(t-\eta)^\alpha}.$$

For the time FFPE, lots of numerical works have been done and most of them concern the case $f(x)=0$ (see e.g., \cite{Jiang2010,Karaa2016,Deng2007,Xu2009,Xu2007,McLean2009,McLean2011,Turner08c,Turner09b}).
Some works on the numerical solution of (\ref{modell}) with $f(x)\neq 0$ are listed as follows: Saadatmandi et al. \cite{Saad1} study a collocation method based on shifted Legendre polynomials
in time and Sinc functions in space;
Fairweather et al. \cite{Fair1}, using the $L_1$ approximation to the time fractional derivative, investigate  an orthogonal
spline collocation method in space;
Deng \cite{Deng} and Cao et al. \cite{Cao1} study numerical schemes for ODE systems derived from spatial semi discretization of (\ref{modell}).

To numerically solve (\ref{modell}), Chen et al. \cite{Chen} present three implicit finite difference schemes and obtain the stability and convergence under conditions that $f(x)$ is a constant or a monotone function. Jiang \cite{Jiang2015} establishes monotone properties
of the numerical schemes given in \cite{Chen}, based on which, a new
proof of stability and convergence is provided under a weaker condition $f(x)\in C[a,b]$. Recently, Vong and Wang \cite{Vong1} develop a
high order compact difference scheme for (\ref{modell}), and obtain its stability and convergence.
In both \cite{Jiang2015} and \cite{Vong1}, the stability and convergence are derived under the condition $h|f(x)|\leq C$ ($h$ denotes space step size and $C$ is a constant),
which requires sufficiently fine space grid in order to obtain desirable numerical solutions when $|f(x)|$ is large.

Finite volume (FV) methods, after being successfully applied to solving integer partial differential equations (see e.g. \cite{Versteeg}),
have been concerned for solving fractional equations in some literatures: for example, \cite{Feng2015,Hejazi2014,Liu2014,Yang2014} study FV methods for space fractional equations; \cite{Hejazi2013} studies FV methods for time-space fractional equations; \cite{Karaa2016} investigates a FV element method for a two dimensional time FPDE (the two dimensional case of (\ref{modell}) with $f=0$).

We develop a FV method for solving (\ref{modell}), using center and upwind differences for space discretization and $L_1$ approximation to the time fractional derivative. Under a discrete $L^1$-norm, we theoretically prove that our method is unconditionally stable and its convergence order is $O(h+\Delta t^{2-\alpha})$ (the convergence order can be improved to $O(h^2+\Delta t^{2-\alpha})$ if the space grid becomes sufficiently fine). Numerical results are given to support our theoretical findings.
  One characteristic of our method is that it has monotone property such that it keeps the
nonnegativity of some physical variables such as density, concentration, etc. and  is efficient even on coarse grids for the case that $|f(x)|$ is large.

In this paper,  $C$ denotes a generic positive constant independent of grid sizes, which may take on different values at different
places.
For easy statement,  we assume that the solution $w$ is sufficiently smooth and $f(x)$ satisfies the following Lipschitz condition
\begin{equation}\label{Lipschitz condition}
|f(x)-f(x')|\leq C|x-x'|, \quad x,x'\in [a,b].
\end{equation}

The rest of this paper is organized as follows: in section 2, we introduce a finite volume scheme for solving (\ref{modell}); in section 3, we show its monotone property, based on which we prove the stability and convergence of the scheme; in section 4, we carry out some numerical tests to support our theory.

\section{Discretisation}
For a positive integer $N$, define the space step size $h=(b-a)/(N+1)$ and uniformly-spaced nodes $x_i=a+ih,i=0,1,\ldots,N+1$. Then we have $N$ control volumes $[x_{i-\frac12},x_{i+\frac12}]$, $i=1,2,\ldots,N$. For a positive integer $L$, define $\Delta t=T/L$. We use the uniform time grid $t_k=k\Delta t (k=0,1,\ldots,L)$ for discretization of the time fractional derivative.

For simplicity, we use the following conventions: for function $g(x,t)$, $g_{i+\frac12,n}=(g)_{(x_{i+\frac12},t_n)}:=g(x_{i+\frac12},t_n)$; for function $g(x)$, $g_{i+\frac12}=(g)_{x_{i+\frac12}}:=g(x_{i+\frac12})$.  Decompose the function $f(x)$ as
$$f=f^u+f^l+f^m,$$
with $f^u(x):=\max(f(x)-\frac{2k_\alpha}h,0)$, $f^l(x):=\min(f(x)+\frac{2k_\alpha}h,0)$, and $f^m:=f-f^u-f^l$ for $x\in [a,b]$.
It is not hard to see that $f^u,f^l,$ and $f^m$ each satisfy the Lipschitz condition with the same positive $C$ as that in (\ref{Lipschitz condition}) and
\begin{equation}\label{eqn add1}
f^u(x)\geq 0,\quad f^l(x)\leq 0,  \quad    -\frac{2k_\alpha}{h}\leq f^m(x)\leq\frac{2k_\alpha}{h}.
\end{equation}

\begin{remark}\label{remark add1}
In designing the numerical method that follows, we only use the values of functions $f^m$, $f^u$ and $f^l$ at points $x_{i+\frac12}$ ($i=0,1,\ldots,N$), i.e., $f^m_{i+\frac12}=f^m(x_{i+\frac12}),f^u_{i+\frac12}=f^u(x_{i+\frac12}),f^l_{i+\frac12}=f^l(x_{i+\frac12})$.
\end{remark}

Evaluating (\ref{modell}) at $t=t_n$ $(n=1,2,\ldots,L)$ and integrating equation (\ref{modell}) over the $i$th control volume $[x_{i-\frac12},x_{i+\frac12}]$ give
\begin{eqnarray}\label{eqn conserve}
\int_{x_{i-\frac12}}^{x_{i+\frac12}}\left.\frac{\partial^\alpha w}{\partial t^\alpha}\right|_{t_n} dx&=&k_\alpha\left(\left(\frac{\partial w}{\partial x}\right)_{{i+\frac12},n} -\left(\frac{\partial w}{\partial x}\right)_{{i-\frac12},n} \right)\nonumber\\
&&\quad-\quad\left((fw)_{{i+\frac12},n}-(fw)_{i-\frac12,n}\right),
\end{eqnarray}
$i=1,2,\ldots,N$.
Denote $w(x_i,t_n)$ by $w_i^n$ ($i=0,\ldots,N+1,n=0,1,\ldots,L$). The left hand term of (\ref{eqn conserve}) may be written as
\begin{eqnarray}\label{eqn time approximation}
&&\int_{x_{i-\frac12}}^{x_{i+\frac12}}\left.\frac{\partial^\alpha w}{\partial t^\alpha}\right|_{t_n} dx=\left(\frac{\partial^\alpha w}{\partial t^\alpha}\right)_{(x_i,t_n)}h-h\gamma^{(1)}_{i,n}\nonumber\\
&=&\frac{h\Delta t^{-\alpha}}{\Gamma(2-\alpha)}\left(w_i^n-\sum\limits_{k=1}^{n-1}(a_{n-k-1}-a_{n-k})w_i^k-a_{n-1}w_i^0\right)-h\gamma^{(1)}_{i,n}-h\gamma^{(2)}_{i,n},
\end{eqnarray}
where the mid rectangular integral formula is used in the first equality, the $L_1$ approximation is used in the second equality
with $a_k=(k+1)^{1-\alpha}-k^{1-\alpha}$ ($k=0,1,2,\ldots$)(see \cite{Oldham,Langlands2005}), and the truncation
errors
$$\gamma^{(1)}_{i,n}:=\left(\left(\frac{\partial^\alpha w}{\partial t^\alpha}\right)_{(x_i,t_n)}h-\int_{x_{i-\frac12}}^{x_{i+\frac12}}\left.\frac{\partial^\alpha w}{\partial t^\alpha}\right|_{t_n} dx\right)/h,$$
$$\gamma^{(2)}_{i,n}:=\left(\frac{h\Delta t^{-\alpha}}{\Gamma(2-\alpha)}\left(w_i^n-\sum\limits_{k=1}^{n-1}(a_{n-k-1}-a_{n-k})w_i^k-a_{n-1}w_i^0\right)-h\left(\frac{\partial^\alpha w}{\partial t^\alpha}\right)_{(x_i,t_n)}\right)/h.$$
We have the truncation errors bounds
\begin{equation}\label{truncation error 12}
h\left|\gamma^{(1)}_{i,n}\right|\leq Ch^3,\quad h\left|\gamma^{(2)}_{i,n}\right|\leq Ch\Delta t^{2-\alpha},\quad i=1,\ldots,N;n=1,\ldots,L,
\end{equation}
where for the second error bound please refer to \cite{Langlands2005}.
The first term on the right hand side of (\ref{eqn conserve}) can be written as
\begin{eqnarray}\label{eqn diffusion 11}
&&k_\alpha\left(\left(\frac{\partial w}{\partial x}\right)_{{i+\frac12},n} -\left(\frac{\partial w}{\partial x}\right)_{{i-\frac12},n}\right)\nonumber\\
&=&k_\alpha\left(\widehat{\left(\frac{\partial w}{\partial x}\right)}_{{i+\frac12},n} -\widehat{\left(\frac{\partial w}{\partial x}\right)}_{{i-\frac12},n}\right)+h\gamma^{(3)}_{i,n},
\end{eqnarray}
where
\begin{eqnarray}\label{eqn diffusion}
\widehat{\left(\frac{\partial w}{\partial x}\right)}_{{i+\frac12},n}:=\frac{w^n_{i+1}-w^n_i}{h},\quad i=0,2,\ldots,N;
\end{eqnarray}
 we have used the midpoint difference formulas for approximations of the derivatives; By Taylor formula, it is easy to see that
\begin{equation}\label{eqn truncation}
h\left|\gamma^{(3)}_{i,n}\right|\leq Ch^3,\quad i=1,2,\ldots,N.
\end{equation}

The second term on the right hand side of (\ref{eqn conserve}) is related to the convection speed $f(x)$. We combine the central and upwind differences  to rewrite
\begin{eqnarray}\label{eqn convection}
&&(fw)_{{i+\frac12},n}=(f^mw)_{{i+\frac12},n}+(f^uw)_{{i+\frac12},n}+(f^lw)_{{i+\frac12},n}\nonumber\\
&=&\widetilde{(f^mw)}_{i+\frac12,n}+\widetilde{(f^uw)}_{i+\frac12,n}+\widetilde{(f^lw)}_{i+\frac12,n}+r^m_{i+\frac12,n}+r^u_{i+\frac12,n}+r^l_{i+\frac12,n}
\end{eqnarray}
for $i=0,1,\ldots, N$, with
\begin{equation}\label{eqn convection11}
\widetilde{(f^mw)}_{i+\frac12,n}:=f^m_{i+\frac12}\frac{w^n_{i}+w^n_{i+1}}2,
\end{equation}
\begin{equation}\label{eqn convection12}
\widetilde{(f^uw)}_{i+\frac12,n}:=f^u_{i+\frac12}w^n_{i},\quad \widetilde{(f^lw)}_{i+\frac12,n}:=f^l_{i+\frac12}w^n_{i+1},
\end{equation}
\begin{equation}\label{eqn convection21}
r^m_{i+\frac12,n}:=(f^mw)_{i+\frac12,n}-\widetilde{(f^mw)}_{i+\frac12,n}=-f^m_{i+\frac12}\left(\frac{\partial^2 w}{\partial x^2}\right)_{i+\frac12,n}\frac{h^2}8+O(h^4),
\end{equation}
\begin{equation}\label{eqn convection22}
r^u_{i+\frac12,n}:=(f^uw)_{i+\frac12,n}-\widetilde{(f^uw)}_{i+\frac12,n}=f^u_{i+\frac12}\left(\frac{\partial w}{\partial x}\right)_{i+\frac12,n}\frac{h}2+O(h^2),
\end{equation}
\begin{equation}\label{eqn convection23}
r^l_{i+\frac12,n}:=(f^lw)_{i+\frac12,n}-\widetilde{(f^lw)}_{i+\frac12,n}=-f^l_{i+\frac12}\left(\frac{\partial w}{\partial x}\right)_{i+\frac12,n}\frac{h}2+O(h^2).
\end{equation}
where we have used Taylor formula in (\ref{eqn convection21}), (\ref{eqn convection22}) and (\ref{eqn convection23}).
Then the second term on the right hand of (\ref{eqn conserve}) may be rewritten as
\begin{equation}\label{eqn convection term}
(fw)_{{i+\frac12},n}-(fw)_{{i-\frac12},n}=\widetilde{(fw)}_{i+\frac12,n}-\widetilde{(fw)}_{i-\frac12,n}-h\gamma^{(4)}_{i,n},
\end{equation}
with
\begin{equation}\label{eqn add3}
\widetilde{(fw)}_{i+\frac12,n}=\widetilde{(f^mw)}_{i+\frac12,n}+\widetilde{(f^uw)}_{i+\frac12,n}+\widetilde{(f^lw)}_{i+\frac12,n}
\end{equation}
and
\begin{equation}\label{eqn add2}
h\gamma^{(4)}_{i,n}=(r^m_{i-\frac12,n}-r^m_{i+\frac12,n})+(r^u_{i-\frac12,n}-r^u_{i+\frac12,n})+(r^l_{i-\frac12,n}-r^l_{i+\frac12,n}).
\end{equation}
\begin{remark}\label{remark kk} It is easy to deduce that for $i=1,\ldots,N$
$$|r^m_{i-\frac12,n}-r^m_{i+\frac12,n}|\leq Ch^3,$$
and
$$|r^u_{i-\frac12,n}-r^u_{i+\frac12,n}|\leq Ch^2,\quad |r^l_{i-\frac12,n}-r^l_{i+\frac12,n}|\leq Ch^2.$$
By the definitions of $f^u$ and $f^l$, $|r^u_{i-\frac12,n}-r^u_{i+\frac12,n}|$ and $|r^l_{i-\frac12,n}-r^l_{i+\frac12,n}|$ become smaller as $h$ becomes smaller, and if ${|f(x)|}\leq \frac{2k_\alpha} h$, they become zeros since $f^u=f^l=0$.
\end{remark}

Substituting (\ref{eqn convection term}), (\ref{eqn diffusion 11}) and (\ref{eqn time approximation}) into (\ref{eqn conserve}) gives
\begin{eqnarray}\label{eqn discrete p}
&&\frac{h\Delta t^{-\alpha}}{\Gamma(2-\alpha)}\left(w_i^n-\sum\limits_{k=1}^{n-1}(a_{n-k-1}-a_{n-k})w_i^k-a_{n-1}w_i^0\right)\nonumber\\
&=&k_\alpha\left(\widehat{\left(\frac{\partial w}{\partial x}\right)}_{{i+\frac12},n} -\widehat{\left(\frac{\partial w}{\partial x}\right)}_{{i-\frac12},n}\right)-\left(\widetilde{(fw)}_{i+\frac12,n}-\widetilde{(fw)}_{i-\frac12,n}\right)+hr^n_{i},
\end{eqnarray}
 $i=1,2,\ldots,N$, $n=1,2,\ldots,L$, where $r_i^n=\sum_{j=1}^4\gamma^{(j)}_{i,n}$.
We use $W^n_i$ to denote the approximation of $w^n_i$. By (\ref{eqn discrete p}),
 we derive the following finite volume (FV) scheme: for $i=1,2,\ldots,N$, $n=1,2,\ldots,L$,
\begin{eqnarray}\label{eqn discrete e}
&&\frac{h\Delta t^{-\alpha}}{\Gamma(2-\alpha)}\left(W_i^n-\sum\limits_{k=1}^{n-1}(a_{n-k-1}-a_{n-k})W_i^k-a_{n-1}W_i^0\right)\nonumber\\
&&=k_\alpha\left(\widehat{\left(\frac{\partial W}{\partial x}\right)}_{{i+\frac12},n} -\widehat{\left(\frac{\partial W}{\partial x}\right)}_{{i-\frac12},n}\right)-\left(\widetilde{(fW)}_{i+\frac12,n}-\widetilde{(fW)}_{i-\frac12,n}\right),
\end{eqnarray}
subject to boundary values $W^n_0=g_1(t_n)$, $W^n_{N+1}=g_2(t_n)$, $n=1,2,\ldots,L$ and initial values $W^0_i$ are taken as the approximations to $w^0_i=\phi(x_i)$, $i=1,\ldots,N$,
where $\widehat{\left(\frac{\partial W}{\partial x}\right)}_{{i+\frac12},n}$ are defined by directly replacing $w$ with $W$ in (\ref{eqn diffusion}), $\widetilde{(fW)}_{i+\frac12,n}$ are defined by directly replacing $w$ with $W$ in (\ref{eqn add3}), (\ref{eqn convection11}) and (\ref{eqn convection12}).

The matrix form of the FV scheme is
\begin{eqnarray}\label{eqn discrete}
\left(\frac{h\Delta t^{-\alpha}}{\Gamma(2-\alpha)}I+A+B\right)W^n=\frac{ h\Delta t^{-\alpha}}{\Gamma(2-\alpha)}\left(\sum\limits_{k=1}^{n-1}(a_{n-k-1}-a_{n-k})W^k+a_{n-1}W^0\right)+d^n,
\end{eqnarray}
$n=1,2,\ldots,L$, where $W^k=(W^k_1,\ldots,W^k_N)^T$, $d^n=(d^n_1,d^n_2,\ldots,d^n_N)\in R^N$, $I\in R^{N\times N}$ is the identity matrix, $A=(a_{ij})\in R^{N\times N}$ is the coefficient matrix corresponding to the first term of the right hand side (\ref{eqn discrete e}), $B=(b_{ij})\in R^{N\times N}$ is the coefficient matrix corresponding to the second term of the right hand side (\ref{eqn discrete e}), the entries of $A,B,d^n$ are listed as follows:

\noindent the $1$th column of $A$
\begin{equation}\label{eqn coff A c1}
a_{11}=\frac{2k_\alpha}h,\quad a_{21}=-\frac{k_\alpha}h,\quad a_{i1}=0, i\geq 3;
\end{equation}

\noindent the $N$th column of $A$
\begin{equation}\label{eqn coff A c2}
a_{NN}=\frac{2k_\alpha}h,\quad a_{(N-1)N}=-\frac{k_\alpha}h, \quad a_{iN}=0, i\leq N-2;
\end{equation}

\noindent the $j$th column of $A$ ($j=2,\ldots,N-1$)
\begin{equation}\label{eqn  A coff2}
a_{jj}=\frac{2k_\alpha}h, \quad a_{(j+1)j}=a_{(j-1)j}=-\frac{k_\alpha}h, \quad a_{ij}=0, i\geq j+2 \hbox{ or }i\leq j-2;
\end{equation}

\noindent the $1$th column of $B$
\begin{equation}\label{eqn coff B c11}
b_{11}=-\frac{f^m_{\frac12}}2+\frac{f^m_{\frac32}}2
+\underline{(-f^l_{\frac12}})+\underline{f^u_{\frac32}},
\end{equation}
\begin{equation}\label{eqn coff B c12}
b_{21}=-\frac{f^m_{\frac32}}2
-f^u_{\frac32}, \quad b_{i1}=0, i\geq 3;
\end{equation}

\noindent the $N$th column of $B$
\begin{eqnarray}\label{eqn coff B cN1}
 b_{NN}=-\frac{f^m_{N-\frac12}}2+\frac{f^m_{N+\frac12}}2+\underline{(-f^l_{N-\frac12})}+\underline{f^u_{N+\frac12}},
\end{eqnarray}
\begin{equation}\label{eqn coff B cN2}
b_{(N-1)N}=\frac{f^m_{N-\frac12}}2
-(-f^l_{N-\frac12}),\quad b_{iN}=0, i\leq N-2;
\end{equation}

\noindent the $j$th column of $B$ ($j=2,\ldots,N-1$)
\begin{equation}\label{eqn  B coff cj1}
b_{jj}=-\frac{f^m_{j-\frac12}}2+\frac{f^m_{j+\frac12}}2
+\underline{(-f^l_{j-\frac12})}+\underline{f^u_{j+\frac12}},
\end{equation}
\begin{equation}\label{eqn  B coff cj2}
b_{(j+1)j}=-\frac{f^m_{j+\frac12}}2
-f^u_{j+\frac12},
\end{equation}
\begin{equation}\label{eqn  B coff cj3}
b_{(j-1)j}=\frac{f^m_{j-\frac12}}2-(-f^l_{j-\frac12});
\end{equation}

\noindent the entries of $d^n$
\begin{equation}\label{eqn coff d1}
d^n_1=\left(\frac{f^m_{\frac12}}2+f^u_{\frac12}+\frac{k_\alpha}{h}\right)g_1(t_n),\quad d^n_i=0, i=2,\ldots,N-1,
\end{equation}
\begin{equation}\label{eqn coff dN}
d^n_{N}=\left(-\frac{f^m_{N+\frac12}}2+(-f^l_{N+\frac12})+\frac{k_\alpha}{h}\right)g_2(t_n).
\end{equation}

Before proceeding forward, we introduce the definition of $M$-matrix and a relevant lemma (see \cite{Berman,Cottle}).
\begin{definition}\label{def1}
A square matrix $\mathbb{A}$ is called an $M$-matrix if it has non-positive off-diagonals and is non-singular with $\mathbb{A}^{-1}\geq 0$.
\end{definition}
\begin{lemma}\label{lem}
If a square matrix $\mathbb{A}$ with non-positive off-diagonals and positive diagonals is strictly row/column diagonal dominant, then $\mathbb{A}$ is an $M$-matrix.
\end{lemma}

The following Lemma is crucial for analyzing the stability and convergence of our scheme.
\begin{lemma}\label{lemma M matrix}
$\frac{h\Delta t^{-\alpha}}{\Gamma(2-\alpha)}I+A+B$ is an M-matrix.
\end{lemma}
\proof
By (\ref{eqn add1}), the underlined quantities in (\ref{eqn coff B c11}), (\ref{eqn coff B cN1}) and (\ref{eqn  B coff cj1}) are nonnegative and $\left|\frac{f^m_{k+\frac12}}2\right|\leq \frac{k_\alpha}{h}, k=0,1,\ldots,N.$
Then through a careful check, we can easily see that: 1. the diagonal entries of $A+B$ are nonnegative; 2. the off-diagonal entries of $A+B$ are
non-positive; 3. $\sum_{i=1}^N(a_{ij}+b_{ij})=0$ for $2\leq j\leq N-1$ and $\sum_{i=1}^N(a_{ij}+b_{ij})\geq0$ for $j=1,N$. So the matrix $\frac{h\Delta t^{-\alpha}}{\Gamma(2-\alpha)}I+A+B$ is strictly column dominant with non-positive off-diagonals and positive diagonals, and then $\frac{h\Delta t^{-\alpha}}{\Gamma(2-\alpha)}I+A+B$ is an $M$-matrix by Lemma \ref{lem}.\eproof

\begin{remark}\label{remark 4}
When we view $w$ as the concentration of some kind of gas, the considered problem is just a mass transfer model and the solution $w$ is nonnegative.
It is easy to see that $a_{n-k-1}-a_{n-k}$ and $a_{n-1}$ in (\ref{eqn discrete}) are positive, and the coefficients of $g_1(t_n)$ in (\ref{eqn coff d1}) and  $g_2(t_n)$ in (\ref{eqn coff dN}) are nonnegative.
 By induction we may conclude that the FV scheme generates nonnegative solution $W^n$ since in practice the boundary values $W_{0}^n$, $W_{N+1}^n$ ($n=1,2,\ldots,L$) and the initial values $W^0_i$ ($i=1,2,\ldots,N$) are all nonnegative. So our FV scheme can keep
 the nonnegativity of the concentration.
\end{remark}

\section{Stability and convergence}
We investigate the stability and convergence of our FV scheme (\ref{eqn discrete e}), whose matrix form is given in (\ref{eqn discrete}).
Denote $e^n_i=w^n_i-W^n_i$ ($i=0,1,\ldots,N+1;n=0,1,\ldots,L$),  $w^n=(w^n_1,w^n_2,\ldots,w^n_N)^T$ and $e^n=(e^n_1,e^n_2,\ldots,e^n_N)^T$. For $W=(W_1,W_2,\ldots,W_N)^T\in R^N$ define the following discrete $L_1$ norm
$$||W||_1:=\sum\limits_{i=1}^{N}h|W_i|.$$
Subtracting (\ref{eqn discrete e}) from (\ref{eqn discrete p}) gives the error equations of the FV scheme
\begin{eqnarray}\label{eqn error d}
&&\frac{h\Delta t^{-\alpha}}{\Gamma(2-\alpha)}\left(e_i^n-\sum\limits_{k=1}^{n-1}(a_{n-k-1}-a_{n-k})e_i^k-a_{n-1}e_i^0\right)\nonumber\\
&=&k_\alpha\left(\widehat{\left(\frac{\partial e}{\partial x}\right)}_{{i+\frac12},n} -\widehat{\left(\frac{\partial e}{\partial x}\right)}_{{i-\frac12},n}\right)-\left(\widetilde{(fe)}_{i+\frac12,n}-\widetilde{(fe)}_{i-\frac12,n}\right)+hr^n_i,
\end{eqnarray}
$i=1,\ldots,N,n=1,\ldots,L$, with $e^n_0=e^n_{N+1}=0$ ($n=1,2,\ldots,L$), $\widehat{\left(\frac{\partial e}{\partial x}\right)}_{{i+\frac12},n}$ being defined by directly replacing $w$ with $e$ in (\ref{eqn diffusion}), $\widetilde{(fe)}_{i+\frac12,n}$ being defined by directly replacing $w$ with $e$ in (\ref{eqn add3}), (\ref{eqn convection11}) and (\ref{eqn convection12}).

The matrix form of the error equations (\ref{eqn error d}) is
\begin{eqnarray}\label{eqn error M}
\left(\frac{h\Delta t^{-\alpha}}{\Gamma(2-\alpha)}I+A+B\right)e^n=\frac{ h\Delta t^{-\alpha}}{\Gamma(2-\alpha)}\left(\sum\limits_{k=1}^{n-1}(a_{n-k-1}-a_{n-k})e^k+a_{n-1}e^0\right)+hr^n,
\end{eqnarray}
$n=1,2,\ldots,L$, where $r^n=(r^n_1,r_2^n,\ldots,r_N^n)^T$.

In the following, we investigate the monotone property of the FV scheme.
\begin{lemma}\label{lemma2}
We have
  $$e^n\geq 0\quad \hbox{ for }n=1,\ldots,L, $$
  if $e^n=(e^n_1,\ldots,e^n_{N})^T$ ($n=1,\ldots,L$) is generated by (\ref{eqn error M}) with  $e^0=(e^0_1,\ldots,e^0_{N})^T\geq 0$ and $r^n=(r^n_1,\ldots,r^n_{N})^T\geq 0$ ($n=1,\ldots,L$).
\end{lemma}
\proof By Lemma \ref{lemma M matrix}, $\frac{h\Delta t^{-\alpha}}{\Gamma(2-\alpha)}I+A+B$ is an $M$-matrix, and so  $\left(\frac{h\Delta t^{-\alpha}}{\Gamma(2-\alpha)}I+A+B\right)^{-1}\geq 0$.
From
 $$a_{n-1}>0,\quad a_{n-k-1}-a_{n-k}>0,\quad k=1,\ldots,n-1,n=1,\ldots,L,$$
 the lemma can be proved by induction.\eproof

The following corollary is obtained directly from Lemma \ref{lemma2}.
\begin{co}\label{co1}
Let $\bar\epsilon,\tilde\epsilon\in R^{N}$, $\bar r^n,\tilde r^n\in R^{N}$ satisfy $\bar\epsilon\geq\tilde\epsilon$, $\bar r^n\geq \tilde r^n$ $(n=1,\ldots,L)$. We have
$$\bar{e}^n\geq  \tilde{e}^n,\quad \hbox{for }n=1,2,\ldots,L,$$ if $\{\bar{e}^{n}\}_{n=1}^{L}$ and $\{\tilde{e}^{n}\}_{n=1}^{L}$ are generated by (\ref{eqn error M}) with $\{e^0=\bar\epsilon,r^n=\bar{r}_n\}$ and with $\{e^0=\tilde\epsilon,r^n=\tilde r_n\}$, respectively.
\end{co}

In the analysis of stability and convergence, we shall use the above monotone results and the following lemma.
\begin{lemma}\label{lemma convergence}\cite{Jiang2015}
Let $b^{(n)}_i$($i=0,1\ldots,n;n=1,2\ldots,L$) be positives with $b^{(n)}_n\geq\sum_{i=0}^{n-1}b^{(n)}_i$, and positives $\varepsilon_n$, $n=1,\ldots,L$, satisfy
\begin{equation}\label{eqn3}
b_n^{(n)}\varepsilon_n\leq \sum_{k=1}^{n-1}b_k^{(n)}\varepsilon_{k}+b_{0}^{(n)}\mu+\kappa_n,\quad n=1,\ldots,L,
\end{equation}
where $\mu,\kappa_n,n=1,\ldots,L$ are positives. Then we have
\begin{equation}\label{eqn22}
\varepsilon_n\leq \mu+{\max\limits_{1\leq i\leq n}\kappa_i}/{\min\limits_{1\leq i\leq n}b_0^{(i)}},\quad n=1,\ldots,L.
\end{equation}
\end{lemma}

Finally the stability and convergence for the FV scheme are the results of the following lemma.
\begin{lemma}\label{lemma error estimate}
 Given $e^0=(e^0_1,\ldots,e^0_{N})^T,r^n=(r^n_1,\ldots,r^n_N)^T\in R^{N}$, $n=1,\ldots,L$, let $e^n=(e^n_1,\ldots,e^n_{N})^T\in R^{N}$, $n=1,\ldots,L$,  be generated by (\ref{eqn error M}). Then we have
\begin{equation}\label{eqn33}
||e^n||_1\leq ||e^0||_1+Cr, \quad n=1,2,\ldots,L,
\end{equation}
where $r=\max\limits_{1\leq n\leq L}||r^n||_1$.
\end{lemma}
\proof
Here only we prove that
\begin{equation}\label{eqn5}
||e^n||_1\leq\sum\limits_{k=1}^{n-1}(a_{n-k-1}-a_{n-k})||e^{k}||_1+a_{n-1}||e^0||_1+\Gamma(2-\alpha)\Delta t^{\alpha}||r^n||_1
\end{equation}
 if $\{e^0\geq 0,r^n\geq 0\}$, $n=1,2,\ldots,L$. The rest of proof is the same as that for Theorem 3.1 in \cite{Jiang2015}.
Assume that $\{e^0\geq 0,r^n\geq 0\}$, $n=1,2,\ldots,L$. By Lemma \ref{lemma2},
$$e^n\geq 0, \quad n=1,2,\ldots,L.$$
write $e^n_0=e^n_{N+1}=0 (n=1,2,\ldots,L)$. Then we know that $e^n_i(i=0,1,\ldots,N+1;n=1,2,\ldots,L)$ satisfy (\ref{eqn error d}).
Summing (\ref{eqn error d}) over $i=1,2,\ldots,N$ gives
\begin{eqnarray}\label{eqnarray1}
&&\frac{\Delta t^{-\alpha}}{\Gamma(2-\alpha)}\left[||e^n||_1-\sum\limits_{k=1}^{n-1}(a_{n-k-1}-a_{n-k})||e^{k}||_1-a_{n-1}||e^0||_1\right]\nonumber\\
&&=k_\alpha\left(\widehat{\left(\frac{\partial e}{\partial x}\right)}_{{N+\frac12},n} -\widehat{\left(\frac{\partial e}{\partial x}\right)}_{{\frac12},n}\right)-\left(\widetilde{(fe)}_{N+\frac12,n}-\widetilde{(fe)}_{\frac12,n}\right)+||r^n||_1\nonumber\\
&&=\frac{-k_\alpha}{h}e^n_N+\frac{-k_\alpha}{h}e^n_1-\left(\frac{f^m_{N+\frac12}}2+\underline{f^u_{N+\frac12}}\right)e^n_N\nonumber\\
&&\qquad +\left(\frac{f^m_{\frac12}}2-\underline{(-f^l_{\frac12})}\right)e^n_1+||r^n||_1.
\end{eqnarray}
 Noticing that $e_{1}^n\geq 0$, $e_{N-1}^n\geq 0$, $\left|\frac{f^m_{N+\frac12}}2\right|\leq \frac{k_\alpha}h$, $\left|\frac{f^m_{\frac12}}2\right|\leq \frac{k_\alpha}h$ and the underlined quantities are nonnegative, we have
 $$\frac{-k_\alpha}{h}e^n_N+\frac{-k_\alpha}{h}e^n_1-\left(\frac{f^m_{N+\frac12}}2+{f^u_{N+\frac12}}\right)e^n_N
+\left(\frac{f^m_{\frac12}}2-{(-f^l_{\frac12})}\right)e^n_1\leq 0,$$
 and furthermore we obtain (\ref{eqn5}).
\eproof

Now we are in a position to present the main results of this paper.
\begin{theorem}\label{co2}
Let $e^n$ be generated by (\ref{eqn error M}) with any $e^0\in R^N$ and $r^n=0$ $ (n=1,2,\ldots,L)$, and then we have the following unconditional stability for the FV scheme
\begin{equation}\label{eqn333}
||e^n||_1\leq ||e^0||_1, \quad n=1,2,\ldots,L.
\end{equation}
\end{theorem}
\proof The theorem is a direct consequence of Lemma \ref{lemma error estimate}.\eproof

\begin{theorem}\label{co3}  For the FV scheme of solving \{(\ref{modell}),(\ref{p_b1}),(\ref{p_b2})\}, let $e^n_i$, $i=1,\ldots,N$, $n=1,\ldots,L$, be the errors between the real solution $w^n_i$ and the computational solution $W^n_i$ at $(x_i,t_n)$, and $e^n=(e^n_1,\ldots,e^n_{N})^T$. If the solution $w$ is sufficiently smooth, we have, for any $n=1,2,\ldots,L$,
\begin{equation}\label{eqn12}
||e^n||_1\leq ||e^0||_1+C(\tau^{2-\alpha}+h),
\end{equation}
or
\begin{equation}\label{eqn112}
||e^n||_1\leq ||e^0||_1+C(\tau^{2-\alpha}+h^2)
\end{equation}
if the space grid is sufficiently fine such that $h\cdot\max|f(x)|\leq {2k_\alpha}$.
\end{theorem}
\proof By the truncation errors $r_i^n=\sum_{j=1}^4\gamma^{(j)}_{i,n}$, (\ref{truncation error 12}), (\ref{eqn truncation}), (\ref{eqn add2}) and Remark \ref{remark kk}, we know $||r^n||_1\leq C(\tau^{2-\alpha}+h)$ and $||r^n||_1\leq C(\tau^{2-\alpha}+h^2)$ if $h\cdot\max|f(x)|\leq {2k_\alpha}$. Since $\{e^n\}_{n=1}^L$ satisfies (\ref{eqn error M}), the theorem follows directly from Lemma \ref{lemma error estimate}.\eproof

\section{Numerical results}
In this section, we shall numerically test our FV scheme for solving \{(\ref{modell}),(\ref{p_b1}),(\ref{p_b2})\}. Instead of (\ref{modell}), we use the following equation with a source term
\begin{equation}\label{problem test}
\frac{\partial^\alpha w}{\partial t^\alpha}=\left[k_\alpha\frac{\partial^2}{\partial x^2}-\frac{\partial}{\partial x}f(x)\right]w(x,t)+g(x,t), \quad (x,t)\in [0,1]\times[0,1].
\end{equation}
Due to the presence of a new term $g(x,t)$ for the original equation (\ref{modell}), we modify FV scheme (\ref{eqn discrete e}) by adding a term $hg(x_i,t_n)$ ( it is from $\int_{x_{i-\frac12}}^{x_{i+\frac12}}g(x,t_n)dx= hg(x_i,t_n)+O(h^3)$) on its right hand side, which brings no difference to our theoretical analysis.

We here present two examples: one serves to test the convergence rates of our FV scheme, the other to make simple comparisons between our FV Scheme and a finite difference method.
\begin{exam}\label{exm1}
We test convergence rates of our FV scheme for solving \{(\ref{problem test}),(\ref{p_b1}),(\ref{p_b2})\} with $k_\alpha=1$, $f(x)=(x-x^2)+400$, $\phi(x)=0$, $g_1(t)=t^2$, $g_2(t)=-t^2$ and
\begin{eqnarray}\label{eqn ab}
 g(x,t)&=&\frac{\Gamma(3)}{\Gamma(3-\alpha)}t^{2-\alpha}\cos(\pi x)+k_\alpha \pi^2 t^2 \sin(\pi x)\nonumber\\
 &&\qquad +t^2\left[(1-2x)\cos(\pi x)-(x-x^2+400)\pi \sin(\pi x)\right].
\end{eqnarray}
 The exact solution of the problem is given by
\[
 w(x,t)=t^2\cos(\pi x).
\]
\end{exam}

The numerical results testing the convergence rate for space are listed  in Tables~\ref{t1}-\ref{t3}, and those for time are listed in Tables~\ref{t4}-\ref{t6}, where

\[
 \hbox{Rate for space}=\left|\frac{\ln (\|\hbox{Error on finer grid}\|_1/\|\hbox{Error on coarser grid}\|_1)}{\ln (N+1 \hbox{ of finer grid}/N+1 \hbox{ of coarser grid})}\right|,
\]
\[
 \hbox{Rate for time}=\left|\frac{\ln (\|\hbox{Error on finer grid}\|_1/\|\hbox{Error on coarser grid}\|_1)}{\ln (L \hbox{ of finer grid}/L \hbox{ of coarser grid})}\right|.
\]
 The tables show that the convergence rate for time is order $2-\alpha$  and the convergence rate for space is order $1$ which increases to order $2$ when the space grid become sufficiently fine.
  The rate increasing in Tables~\ref{t1}-\ref{t3} when $N$ changes from $N=20$ to $320$ is because the truncation $h\gamma^{(4)}_{i,n}$ changes from being $O(h^2)$ to being $O(h^3)$ continuously as $h$ becomes smaller (see Remark \ref{remark kk}).

\begin{table}[htbp]
\caption{Convergence rate for space with $\alpha=0.2$ and $L=10000$.}
\smallskip
  \centering
  \begin{tabular}{c|c|c|c|c}
    \hline
    $N+1$ & 10 & 20 & 40 & 80 \\
    \hline
    $\max\limits_{n}\|e^{(n)}\|_{\infty}$ & $1.421\times10^{-1}$    &  $6.920\times10^{-2}$  &  $3.119\times10^{-2}$     &    $1.180\times10^{-2}$ \\
    $\max\limits_{n}\|e^{(n)}\|_1$        & $8.783\times10^{-2}$    &  $4.349\times10^{-2}$  &  $1.976\times10^{-2}$     & $7.516\times10^{-3}$     \\
     Conv. rate                           &                         &  1.014                 &     1.138                 & 1.395                   \\
    \hline
     $N+1$  & 160 & 320 & 640 & 1280 \\
    \hline
    $\max\limits_{n}\|e^{(n)}\|_{\infty}$ & $2.012\times10^{-3}$    & $3.206\times10^{-5}$ & $8.012\times10^{-6}$ & $2.003\times10^{-6}$ \\
    $\max\limits_{n}\|e^{(n)}\|_1$        & $1.295\times10^{-3}$    & $1.600\times10^{-5}$ & $4.001\times10^{-6}$ & $1.000\times10^{-6}$  \\
     Conv. rate                           & 2.537                   & 6.339                & 2.000               & 2.000                \\
    \hline
  \end{tabular}\label{t1}
\end{table}
\begin{table}[htbp]
\caption{Convergence rate for space with $\alpha=0.5$ and $L=10000$.}
\smallskip
  \centering
  \begin{tabular}{c|c|c|c|c}
    \hline
    $N+1$ & 10 & 20 & 40 & 80 \\
    \hline
    $\max\limits_{n}\|e^{(n)}\|_{\infty}$ & $1.421\times10^{-1}$    &  $6.918\times10^{-2}$  &  $3.118\times10^{-2}$     &    $1.180\times10^{-2}$ \\
    $\max\limits_{n}\|e^{(n)}\|_1$        & $8.779\times10^{-2}$    &  $4.347\times10^{-2}$  &  $1.976\times10^{-2}$     & $7.513\times10^{-3}$     \\
     Conv. rate                           &                         &  1.014                 &     1.138                 & 1.395                   \\
    \hline
     $N+1$  & 160 & 320 & 640 & 1280 \\
    \hline
    $\max\limits_{n}\|e^{(n)}\|_{\infty}$ & $2.012\times10^{-3}$    & $3.204\times10^{-5}$ & $8.009\times10^{-6}$ & $2.002\times10^{-6}$ \\
    $\max\limits_{n}\|e^{(n)}\|_1$        & $1.295\times10^{-3}$    & $1.600\times10^{-5}$ & $4.000\times10^{-6}$ & $9.998\times10^{-7}$  \\
     Conv. rate                           & 2.537                   & 6.339                & 2.000               & 2.000                \\
    \hline
  \end{tabular}\label{t2}
\end{table}

\begin{table}[htbp]
\caption{Convergence rate for space with $\alpha=0.8$ and $L=10000$.}
\smallskip
  \centering
  \begin{tabular}{c|c|c|c|c}
    \hline
    $N+1$ & 10 & 20 & 40 & 80 \\
    \hline
    $\max\limits_{n}\|e^{(n)}\|_{\infty}$ & $1.421\times10^{-1}$    &  $6.916\times10^{-2}$  &  $3.117\times10^{-2}$     & $1.180\times10^{-2}$    \\
    $\max\limits_{n}\|e^{(n)}\|_1$        & $8.776\times10^{-2}$    &  $4.345\times10^{-2}$  &  $1.975\times10^{-2}$     & $7.511\times10^{-3}$    \\
     Conv. rate                           &                         &  1.014                 &     1.138                 & 1.395                   \\
    \hline
     $N+1$  & 160 & 320 & 640 & 1280 \\
    \hline
    $\max\limits_{n}\|e^{(n)}\|_{\infty}$ & $2.011\times10^{-3}$    & $3.203\times10^{-5}$ & $8.006\times10^{-6}$ & $2.001\times10^{-6}$ \\
    $\max\limits_{n}\|e^{(n)}\|_1$        & $1.294\times10^{-3}$    & $1.599\times10^{-5}$ & $3.993\times10^{-6}$ & $9.938\times10^{-7}$  \\
     Conv. rate                           & 2.537                   & 6.339                & 2.002                & 2.007                \\
    \hline
  \end{tabular}\label{t3}
\end{table}

\begin{table}[htbp]
\caption{Convergence rate for time with $\alpha=0.2$ and $N=15000$.}
\smallskip
  \centering
  \begin{tabular}{c|c|c|c|c|c}
    \hline
    $L$   & 10 & 20 & 40 & 80 & 160 \\
    \hline
    $\max\limits_{n}\|e^{(n)}\|_{\infty}$   & $2.14\times10^{-6}$ & $6.51\times10^{-7}$ & $1.92\times10^{-7}$ & $5.55\times10^{-8}$  & $1.60\times10^{-7}$\\
    \hline
    $\max\limits_{n}\|e^{(n)}\|_1$          & $1.35\times10^{-6}$ & $4.10\times10^{-7}$ & $1.20\times10^{-7}$ & $3.46\times10^{-8}$  & $9.98\times10^{-8}$\\
    \hline
    Conv. rate   &    & 1.722 & 1.773 & 1.793 & 1.792 \\
    \hline
  \end{tabular}\label{t4}
\end{table}
\begin{table}[htbp]
\caption{Convergence rate for time with $\alpha=0.5$ and $N=5000$.}
\smallskip
  \centering
  \begin{tabular}{c|c|c|c|c|c}
    \hline
    $L$  &  10 & 20 & 40 & 80 & 160 \\
    \hline
    $\max\limits_{n}\|e^{(n)}\|_{\infty}$   & $1.09\times10^{-5}$ & $3.89\times10^{-6}$ & $1.38\times10^{-8}$ & $4.90\times10^{-7}$ & $1.74\times10^{-7}$\\
    \hline
    $\max\limits_{n}\|e^{(n)}\|_1$          & $6.88\times10^{-6}$ & $2.45\times10^{-6}$ & $8.70\times10^{-7}$ & $3.09\times10^{-7}$  & $1.09\times10^{-7}$\\
    \hline
    Conv. rate   &    & 1.490 & 1.494 & 1.494 & 1.495 \\
    \hline
  \end{tabular}\label{t5}
\end{table}

\begin{table}[htbp]
\caption{Convergence rate for time with $\alpha=0.8$ and $N=5000$.}
\smallskip
  \centering
  \begin{tabular}{c|c|c|c|c|c}
    \hline
    $L$  & 10 & 20 & 40 & 80 & 160 \\
    \hline
    $\max\limits_{n}\|e^{(n)}\|_{\infty}$   & $3.76\times10^{-3}$ & $1.64\times10^{-5}$ & $7.13\times10^{-6}$ & $3.11\times10^{-6}$ & $1.35\times10^{-6}$\\
    \hline
    $\max\limits_{n}\|e^{(n)}\|_1$          & $2.38\times10^{-5}$ & $1.04\times10^{-5}$ & $4.52\times10^{-6}$ & $1.97\times10^{-6}$  & $8.57\times10^{-7}$\\
    \hline
    Conv. rate   &    & 1.198 & 1.198 & 1.199 & 1.199 \\
    \hline
  \end{tabular}\label{t6}
\end{table}

\begin{figure}
  \centering
  \includegraphics[width=6cm]{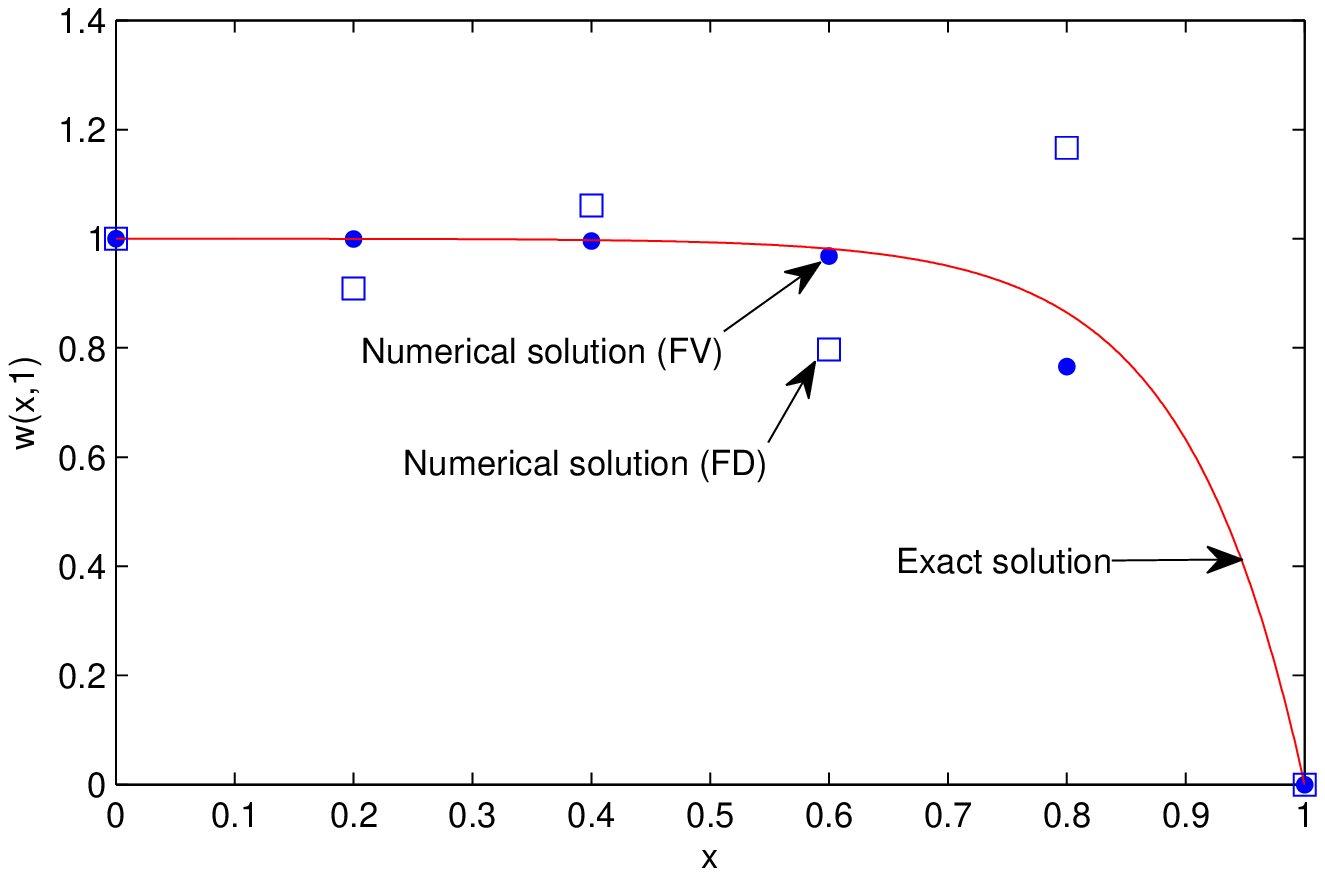}\includegraphics[width=6cm]{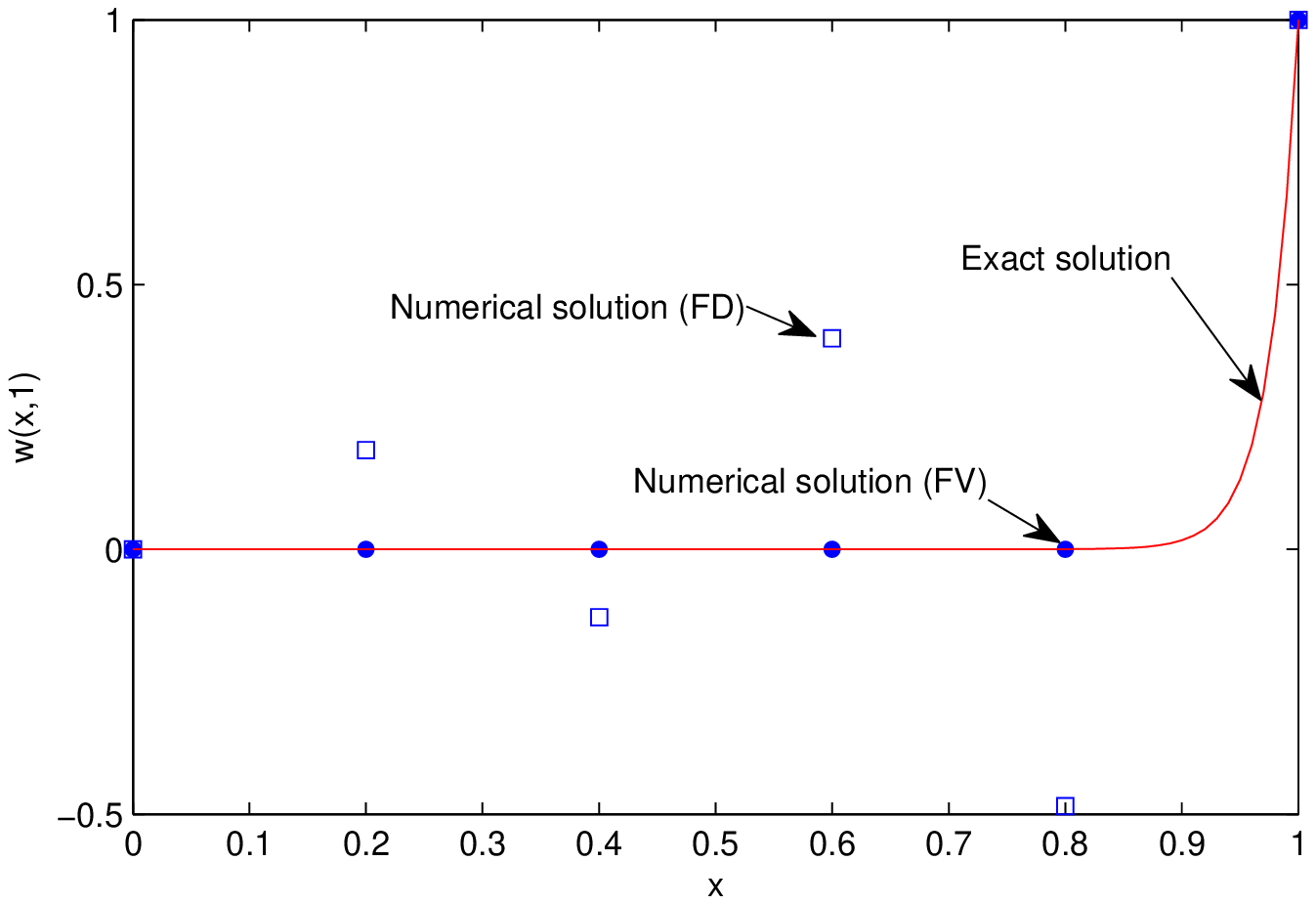}\\
  \caption{Comparisons between the FV scheme and the FD scheme, the left is for Case 1 and the right for Case 2}\label{fig1}
\end{figure}

  As a comparison, we introduce a finite difference (FD) scheme for solving \{(\ref{problem test}),(\ref{p_b1}),(\ref{p_b2})\} which is obtained in the following process: evaluating (\ref{problem test}) at $(x,t)=(x_i,t_n)$, and then the time fractional derivative is discretized in the same way as that in our FV scheme, the space derivatives are discretized using center difference methods which are based on the approximations
$$
\left.\frac{\partial^2 w}{\partial x^2}\right|_{(x_i,t_n)}\approx \frac{w_{i+1}^n-2 w_{i}^n+w_{i-1}^n}{h^2},\quad \left.\frac{\partial (fw)}{\partial x}\right|_{(x_i,t_n)}\approx \frac{f(x_{i+1})w_{i+1}^n-f(x_{i-1})w_{i-1}^n}{2h}.
$$
When the space step size is sufficiently small, the FD scheme has the same monotone property (see \cite{Jiang2015}) as our FV method, and the two schemes have almost the same solution accuracies.
So in the following example, we test the two schemes on relatively coarse space grids, which is meaningful for problems on large space domains.
\begin{exam}\label{exm2}
   We make simple comparisons between our FV scheme and the FD scheme when they are used to solve \{(\ref{problem test}),(\ref{p_b1}),(\ref{p_b2})\}
in the following two cases:

   \noindent Case 1:$\alpha=0.5$, $k_\alpha=1$, $f(x)=(x-x^2)+40$, $\phi(x)=0$, $g_1(t)=t^2$, $g_2(t)=0$ and
\begin{eqnarray}\label{eqn abab1}
 g(x,t)&=&\frac{\Gamma(3)}{\Gamma(3-\alpha)}t^{2-\alpha}\left(1+\frac{1-e^{10x}}{2.20255\times10^4}\right)+\frac{k_\alpha t^2 e^{10x}}{2.20255\times10^2}\nonumber\\
 &&+t^2\left[(1-2x)\left(1+\frac{1-e^{10x}}{2.20255\times10^4}\right)-(x-x^2+400)\frac{e^{10x}}{2.20255\times10^3}\right];
\end{eqnarray}
  \noindent Case 2:$\alpha=0.5$, $k_\alpha=1$, $f(x)=(x-x^2)+40$, $\phi(x)=0$, $g_1(t)=g_2(t)=0$ and
$g(x,t)=0.$

\noindent In Case 1, the exact solution is $u=1+\frac{1-e^{10x}}{2.20255\times10^4}$, and in Case 2, the numerical solution with sufficiently large $N$ and $L$ is used as the exact solution.
\end{exam}
The numerical solutions $w(x,1)$, with $N=4,L=200$, are drawn in Figure \ref{fig1}. From Case 1, we can see that the solution produced by the FD scheme has oscillations. From Case 2, we see that our FV scheme produces nonnegative solution but the FD scheme does not. Example \ref{exm2} shows that our FV scheme can keep the
nonnegativity of physical variables and it has some advantages when used to solve problems on large space domains.

\end{document}